\documentclass[12pt]{amsart}
\usepackage{amsmath,amscd,amssymb,amsfonts}
\def\ver{Sept. 21, 2008, v.5}
\def\ssbull{\raise.2ex\hbox{${\scriptstyle\bullet}$}}
\def\simto{\buildrel\sim\over\longrightarrow}
\def\msum{\hbox{$\sum$}}
\def\mprod{\hbox{$\prod$}}
\def\moplus{\hbox{$\bigoplus$}}
\def\mcup{\hbox{$\bigcup$}}

\def\bs{{\mathbf s}}
\def\bC{{\mathbf C}}

\def\bN{{\mathbf N}}
\def\bP{{\mathbf P}}
\def\bQ{{\mathbf Q}}
\def\bR{{\mathbf R}}
\def\bZ{{\mathbf Z}}
\def\cA{{\mathcal A}}
\def\cD{{\mathcal D}}
\def\cE{{\mathcal E}}
\def\cG{{\mathcal G}}
\def\cH{{\mathcal H}}
\def\cI{{\mathcal I}}
\def\cJ{{\mathcal J}}
\def\cL{{\mathcal L}}
\def\cM{{\mathcal M}}
\def\ocA{\bar{\mathcal A}}
\def\cO{{\mathcal O}}
\def\fa{{\mathfrak a}}
\def\fm{{\mathfrak m}}
\def\tR{\widetilde{R}}
\def\tb{\widetilde{b}}

\def\tF{\widetilde{F}}
\def\tm{\widetilde{m}}
\def\tn{\widetilde{n}}
\def\tD{\widetilde{D}}
\def\tE{\widetilde{E}}
\def\tM{\widetilde{M}}
\def\tP{\widetilde{P}}
\def\tH{\widetilde{H}}
\def\tU{\widetilde{U}}

\def\tX{\widetilde{X}}

\def\tal{\widetilde{\alpha}}
\def\codim{\hbox{\rm codim}}
\def\JN{\hbox{\rm JN}}
\def\Gr{\hbox{\rm Gr}}
\def\DR{\hbox{\rm DR}}

\def\Ker{\hbox{\rm Ker}}

\def\Spec{\hbox{Spec}\,}
\def\cSpec{\hbox{${\mathcal S}pec$}}
\def\Sp{\hbox{\rm Sp}}
\def\hSp{\hbox{\rm $\widehat{\rm S}$p}}

\def\Sing{\hbox{\rm Sing}}
\def\Supp{\hbox{\rm Supp}}
\def\CV{\hbox{\rm CV}}
\def\wdg{\hbox{$\wedge$}}
\def\rd{\partial}
\def\finv{\hbox{$\frac{1}{f}$}}
\def\Mustata{Musta\c{t}\v{a}}
\def\red{{\rm red}}
\def\({{\rm (}}
\def\){{\rm )}}

\begin{document}
\title{On $b$-function, spectrum and multiplier ideals}
\author{Morihiko Saito}
\address{RIMS Kyoto University, Kyoto 606-8502 Japan}
\email{msaito@kurims.kyoto-u.ac.jp}
\date{\ver}
\begin{abstract}
We survey some recent developments in the theory of $b$-function,
spectrum, and multiplier ideals together with certain interesting
relations among them including the case of arbitrary subvarieties.
\end{abstract}
\maketitle

\centerline{Dedicated to Professor Masaki Kashiwara}

\bigskip
\centerline{\bf Introduction}

\bigskip\noindent
It has been known that there are certain interesting relations
between $b$-function, spectrum and multiplier ideals.
We give a survey on this topic.
We first consider the case of hypersurfaces and then
arbitrary subvarieties.
We recall the definition of $b$-function, spectrum and
multiplier ideals, and explain certain properties
together with interesting relations among them.
We also explain the cases of hyperplane arrangements and
monomial ideals.

In Section~1 we recall the definition of $b$-function in the
hypersurface case and explain some related topics including
the $V$-filtration of Kashiwara and Malgrange.
In Section~2 we recall the definition of spectrum in the
hypersurface case and explain some known results mainly due
to Steenbrink.
In Section~3 we recall the definition of multiplier ideals
in the general case and give an extension theorem generalizing
\Mustata's formula in the case of hyperplane arrangements.
In Section~4 we explain certain relations among $b$-function,
spectrum and multiplier ideals in the hypersurface case.
In Section~5 we treat the case of hyperplane arrangements.
In Section~6 we define the $b$-function in the general
case and explain a relation with the multiplier ideals.
In Section~7 we define the spectrum in the general case
and explain a relation with the multiplier ideals.
In Section~8 we treat the monomial ideal case.

\medskip
In this paper we use the following

\medskip\noindent
{\bf Notation.}
$b_f(s)$ = $b$-function of $f$,
$R_f=$ roots of $b_f(-s)$, $\alpha_f=\min R_f$,
$m_{f,\alpha}$ = multiplicity of $\alpha \in R_f$.
Similarly for $\tR_f,\tm_{f,\alpha},\tal_f$ with $b_f(s)$ replaced
by $\tb_f(s):=b_f(s)/(s+1)$.
$R_{f,x},m_{f,x,\alpha},\alpha_{f,x}$ are associated to the local
$b$-function $b_{f,x}(s)$.
$R'_{f,x}=\cup_{y\ne x}R_{f,y}$,
$\alpha'_{f,x}=\min_{y\ne x}\{\alpha_{f,y}\}$.

\bigskip\bigskip
\centerline{\bf 1. $b$-function of a hypersurface}

\bigskip\noindent
In this section we recall the definition of $b$-function in the
hypersurface case and explain some related topics including
the $V$-filtration of Kashiwara and Malgrange.

\medskip\noindent
{\bf 1.1.~ Definition.}
Let $X$ be a complex manifold or a smooth complex algebraic variety,
and $f$ be a holomorphic or algebraic function on $X$.
Let $\cD_X$ be the sheaf of linear differential operators on $X$.
Set $\rd_i=\rd/\rd x_i$ for local coordinates $x_1,\dots,x_n$.
Then
$$
\cD_X[s]f^s\subset \cO_X[{\finv}][s]f^s\quad\text{with}\,\,\,
\rd_i f^s=s(\rd_if)f^{s-1}.
$$
The $b$-function (i.e. the Bernstein-Sato polynomial) $b_f(s)$
is the monic polynomial of the smallest degree such that
$$
b_f(s)f^s=P(x,\rd_x,s)f^{s+1}\quad\text{in}\,\,\,
\cO_X[{\finv}][s]f^s,
$$
where $P(x,\rd_x,s)\in\cD_X[s]$. Locally, this coincides with the
minimal polynomial of the action of $s$ on
$$
\cD_X[s]f^s/\cD_X[s]f^{s+1}.
$$
The latter definition is valid in a more general case.

We define $b_{f,x}(s)$ replacing $\cD_X$ with $\cD_{X,x}$.

\medskip\noindent
{\bf 1.2.~Remark.}
The $b$-function or Bernstein-Sato polynomial for a
hypersurface was introduced by Sato [41] and Bernstein [3],
see also [4].

\medskip\noindent
{\bf 1.3.~Observation.}
Let $i_f:X\to \tX:=X\times\bC$ denote the graph embedding.
Set
$$
\tM=i_{f+}\cO_X=\cD_{\tX}\delta(f-t)=\cO_{X\times\bC}
\bigl[\hbox{$\frac{1}{f-t}$}\bigr]/\cO_{X\times\bC},
\leqno(1.3.1)
$$
This is a free $\cO_X[\rd_t]$-module of rank 1 with basis
$\delta(f-t)$ which is identified with the class of $\frac{1}{f-t}$.
Here $i_{f+}$ denotes the direct image as a $\cD$-module, and
$t$ is the coordinate of $\bC$.
The action of $\rd_i$, $t$ on $\delta(f-t)$ is given by
$$
\rd_i\delta(f-t)=-(\rd_if)\rd_t\delta(f-t),\quad
t\delta(f-t)=f\delta(f-t).
\leqno(1.3.2)
$$
Then $f^s$ is canonically identified with $\delta(f-t)$ by setting
$s = -\rd_tt$, and there is a canonical isomorphism as
$\cD_X[s]$-modules
$$
\cD_X[s]f^s=\cD_X[s]\delta(f-t).
\leqno(1.3.3)
$$

\medskip\noindent
{\bf 1.4.~V-filtration.}
We say that
$V$ is a filtration of Kashiwara [25] and Malgrange [28] along $f$ if
$V$ is exhaustive, separated, and satisfies the following conditions
for any $\alpha\in\bQ$:

\smallskip\noindent
(i) $V^{\alpha}\tM$ is a coherent $\cD_X[s]$-submodule of $\tM$.

\smallskip\noindent
(ii) $tV^{\alpha}\tM\subset V^{\alpha+1}\tM$ and the equality holds for
$\alpha\gg 0$.

\smallskip\noindent
(iii) $\rd_tV^{\alpha}\tM\subset V^{\alpha-1}\tM$.

\smallskip\noindent
(iv) $\rd_tt-\alpha$ is nilpotent on $\Gr_V^{\alpha}\tM$.

\medskip
(If $V$ exists, it is unique.)

\medskip\noindent
{\bf 1.5.~Relation with the $b$-function.} Assume $X$ is affine or
Stein and relatively compact. Then the multiplicity of a root
$\alpha$ of $b_f(s)$ is given by the degree of the minimal polynomial
of $s-\alpha$ on
$$
\Gr_V^{\alpha}(\cD_X[s]f^s/\cD_X[s]f^{s+1}),
\leqno(1.5.1)
$$
using the isomorphism (1.3.3) where $s=-\rd_tt$.
Note that
$V^{\alpha}\tM$ for $\alpha\in\bQ$ and $\cD_X[s]f^{s+i}$ for
$i\in\bN$ are `lattices' of $\tM$, i.e.
$$
V^{\alpha}\tM\subset\cD_X[s]f^{s+i}\subset V^{\beta}\tM\quad
\text{for}\,\,\alpha\gg i\gg \beta,
\leqno(1.5.2)
$$
and $V^{\alpha}\tM$ is an analogue of the Deligne extension [11]
with eigenvalues in $[\alpha,\alpha+1)$.
This is quite similar to the case of differential equations
of one variable with regular singularities.
The existence of $V$ is equivalent to the existence of
$b_f(s)$ locally.

\medskip\noindent
{\bf 1.6.~Theorem} (Kashiwara [24], [25], Malgrange [28]). {\it
The filtration $V$ exists on $\tM:=i_{f+}M$ for any holonomic
$\cD_X$-module $M$ {\rm(}where $V$ is indexed by $\bC$\).
}

\medskip\noindent
{\bf 1.7.~Remarks.}
(i) There are lots of ways to show this theorem.
Indeed, it is essentially equivalent to the existence of
the $b$-function in a generalized sense. In case $M$ is regular,
one way is to use a resolution of singularities and reduce to
the case where the characteristic variety $\CV(M)$ has normal
crossings.

(ii) A holonomic $\cD$-module $M$ is called quasi-unipotent
if the local monodromies of the local systems
$\cH^j\DR(M)|_{S_i}$ are quasi-unipotent where $\{S_i\}$ is a
suitable Whitney stratification.
This condition is equivalent to the condition that
the filtration $V$ along $f$ is indexed
by $\bQ$ for any locally defined function $f$.
Indeed, the last condition is equivalent to the first condition
since the last condition using $V$ is stable by subquotients so that
we can argue by induction on $\dim\Supp\,M$.

\medskip\noindent
{\bf 1.8.~Relation with vanishing cycle functors.}
Let $\rho:X_t\to D$ be a `good' retraction
where $D=f^{-1}(0)$, and $X_t=f^{-1}(t)$ with $t\ne 0$ sufficiently
near 0.
This is obtained by using an embedded resolution of singularities of
$(X,D)$, since the existence of such a retraction is well known in
the normal crossing case and it is enough to compose it with the
blown-down.
Then there are canonical isomorphisms
$$
\psi_f\bC_X=\bR\rho_*\bC_{X_t},\quad
\varphi_f\bC_X=\psi_f\bC_X/\bC_{D},
\leqno(1.8.1)
$$
where $\psi_f\bC_X,\varphi_f\bC_X$ are nearby and vanishing cycle
sheaves, see [13].

Let $F_x$ denote the Milnor fiber around $x\in D$. Then we have
$$
\aligned
(\cH^j\psi_f\bC_X)_x&=H^j(F_x,\bC),\\
(\cH^j\varphi_f\bC_X)_x&=\tH^j(F_x,\bC).
\endaligned
\leqno(1.8.2)
$$

For a $\cD_X$-module $M$ admitting the V-filtration on
$\tM=i_{f+}M$ indexed by $\bQ$, we define $\cD_X$-modules
$$
\psi_f M=\moplus_{0<\alpha\le 1}\Gr_V^{\alpha}\tM,\quad
\varphi_f M=\moplus_{0\le\alpha<1}\Gr_V^{\alpha}\tM.
\leqno(1.8.3)
$$

\medskip\noindent
{\bf 1.9.~Theorem} (Kashiwara [25], Malgrange [28]). {\it
For any quasi-unipotent regular holonomic $\cD_X$-module $M$, we have
the canonical isomorphisms
$$
\aligned
\DR_X\psi_f(M)&=\psi_f\DR_{X}(M)[-1],
\\
\DR_X\varphi_f(M)&=\varphi_f\DR_{X}(M)[-1],
\endaligned
\leqno(1.9.1)
$$
such that $\exp(-2\pi i\rd_tt)$ on the left-hand side corresponds to
the monodromy $T$ on the right-hand side.
}

\medskip\noindent
{\bf 1.10.~Definition.} Set

\smallskip
$R_f= \{$roots of $b_f(-s)\}$,\,\,
$\alpha_f=\min R_f$,\,\,
$m_{f,\alpha}$ : the multiplicity of $\alpha \in R_f$.

\smallskip\noindent
(Similarly for $R_{f,x}$, etc. for $b_{f,x}(s)$.)

\medskip\noindent
{\bf 1.11.~Theorem} (Kashiwara [23]). {\it We have
$R_f\subset \bQ_{>0}$.}

\medskip
(This is proved by using a resolution of singularities.)

\medskip\noindent
{\bf 1.12.~Theorem} (Kashiwara [25], Malgrange [28]). 
{\it We have}

\smallskip\noindent
(i) {\it $e^{-2\pi iR_f}=\{$the eigenvalues of $T$ on $H^j(F_x,\bC)$
for any ${}\qquad\qquad\qquad\quad x\in D, j\in \bZ\}$.}

\smallskip\noindent
(ii) {\it $m_{f,\alpha}\le\min\{i\mid N^i\psi_{f,\lambda}\bC_X= 0\}$
with $\lambda=e^{-2\pi i\alpha}$.
}

\medskip
Here $\psi_{f,\lambda}=\Ker(T_s-\lambda)\subset \psi_f$ in the
abelian category of perverse sheaves {\rm [2]}, and $N=\log T_u$
with $T=T_sT_u$ the Jordan decomposition.

\medskip\noindent
{\bf 1.13.~Remark.}
This is a corollary of the above Theorem (1.9) of Kashiwara
and Malgrange, and is a generalization of a formula of Malgrange~[27]
in the isolated singularity case, see (4.6).

\medskip\noindent
{\bf 1.14.~Microlocal $b$-function.}
Define $\tR_f,\tm_{f,\alpha},\tal_f$ with $b_f(s)$ replaced
by the
{\it microlocal} (or reduced) $b$-function
$$
\tb_f(s):=b_f(s)/(s+1).
\leqno(1.14.1)
$$
By [38], $\tb_f(s)$ is the monic polynomial of the
smallest degree such that
$$
\tb_f(s)\delta(f-t)=\tP\rd_t^{-1}\delta(f-t),
\leqno(1.14.2)
$$
where $\tP\in \cD_X[s,\rd_t^{-1}]$.

Put $n=\dim X$. Then

\medskip\noindent
{\bf 1.15.~Theorem.} {\it We have}
$$\tR_{f} \subset [\tal_{f},n-\tal_{f}],\quad
\tm_{f,\alpha} \le n-\tal_{f} - \alpha + 1.$$

\medskip
(This follows from the filtered duality for $\varphi_f$, see loc.~cit.)

\medskip\noindent
{\bf 1.16.~Remark.}
If $f$ is weighted-homogeneous with an isolated
singularity at the origin, then we have by an unpublished result of
Kashiwara (mentioned in the end of Introduction of [27])
$$
\tR_f=E_f,\quad
\tm_{f,\alpha}=1\,\,\,\text{for}\,\,\,\alpha\in\tR_f,
\leqno(1.16.1)
$$
where $E_f$ is the set of exponents, see (2.1.2) below.
This assertion also follows from a result of Malgrange in loc.~cit.,
see Th.~(4.6) below.

If $f=\sum_i x_i^2$, then $\tal_f=n/2$ and (1.16.1) follows
from the above Theorem~(1.15).

\newpage
\centerline{\bf 2. Spectrum of a hypersurface}

\bigskip\noindent
In this section we recall the definition of spectrum in the
hypersurface case and explain some known results mainly due
to Steenbrink.

\medskip\noindent
{\bf 2.1.~Spectrum.}
Let $f$ be a function on a complex manifold or a smooth complex
algebraic variety $X$ of dimension $n$.
Let $F_x$ denote the Milnor fiber around $x\in D=f^{-1}(0)$.
Following Steenbrink [45], [47] we define the {\it spectrum}
$$
\aligned
\Sp(f,x)&=\Sp(D,x)=\msum_{\alpha>0}\,n_{f,\alpha}t^{\alpha}\,\,\,\,
\hbox{where}\\
n_{f,\alpha}&=\msum_j\,(-1)^{j-n+1}\dim \Gr_F^p\tH^j(F_x,
\bC)_{\lambda}\,\,\,\,\hbox{with}\\
p&=[n-\alpha],\,\,\lambda=\exp(-2\pi i\alpha).
\endaligned
\leqno(2.1.1)
$$
Here $F$ is the Hodge filtration (see [12], [45]) on
$\tH^j(F_x,\bC)_{\lambda}:=\Ker(T_s-\lambda)$ with $T=T_sT_u$ the
Jordan decomposition.
We define the {\it exponents} by
$$
E_f=\{\alpha\in\bQ\mid n_{f,\alpha}\ne 0\}\subset\bQ_{>0}.
\leqno(2.1.2)
$$

\medskip\noindent
{\bf 2.2.~Isolated singularity case.}
In this case we have by [45] symmetry and positivity
$$
n_{f,\alpha}=n_{f,n-\alpha}\ge 0.
\leqno(2.2.1)
$$
Moreover, for $f,g$ on $X,Y$ we have by Scherk-Steenbrink [43]
and Varchenko [49]
$$
\Sp(f+g,(x,y))=\Sp(f,x)\Sp(g,y),
\leqno(2.2.2)
$$
where the product on the right-hand side is taken in $\bQ[t^{1/e}]$
for some $e\in\bZ_{>0}$.
This can be extended to the non-isolated singularity case
(unpublished).

\medskip\noindent
{\bf 2.3.~Weighted homogeneous isolated singularity case.}
Assume $f$ is weighted homogeneous with positive weights
$w_1,\dots,w_n$, i.e. $f=\sum_{\nu} c_{\nu}x^{\nu}$ with
$c_{\nu}=0$ for $\sum_i w_i \nu_i\ne 1$.
Assume further Sing$\,D=\{0\}$.
Then we have by Steenbrink [44]
$$
\Sp(f,x)=\mprod_i\,(t-t^{w_i})/(t^{w_i}-1).
\leqno(2.3.1)
$$
Indeed, he showed that the left-hand side is given by the
Poincare polynomial of the graded vector space
$$
\Omega_X^n/df\wedge\Omega_X^{n-1},
\leqno(2.3.2)
$$
and it is well known that the latter is calculated by using the
morphism $(f_1,\dots,f_n):\bC^n\to\bC^n$
(where $f_i=\partial f/\partial x_i$).

\medskip\noindent
{\bf 2.4.~Nondegenerate Newton boundary case.}
If $n=2$ and $f$ has nondegenerate Newton boundary $\rd P_f$
such that $\bR^2_{\ge 0}\setminus P_f$ is bounded,
then by Steenbrink [45]
$$
\aligned
&E_f\cap(0,1]=\mcup_{\sigma} E_{\sigma}^{\le 1}\quad\hbox{with}
\\
&E_{\sigma}^{\le 1}=\{L_{\sigma}(u)\mid u\in
\bZ_{>0}^2\cap(\{0\}\cup\sigma)^{\rm conv.hull} \},
\endaligned
\leqno(2.4.1)
$$
were $L_{\sigma}$ is a linear function such that
$L_{\sigma}^{-1}(1)\supset\sigma$.
Here the symmetry of $E_f$ with center 1 is used, see (2.2.1).

For $n>2$, the filtration $V$ on
$\Omega_X^n/df\wedge\Omega_X^{n-1}$ is induced by the Newton
filtration,
and there is a combinatorial description by Steenbrink [45],
see also [33], [51].
(Note that [33] was the origin of the theory of bifiltered
strict complexes.)

\medskip\noindent
{\bf 2.5.~Semicontinuity} (Steenbrink [46]).
For a deformation $\{f_{\lambda}\}_{\lambda\in\Delta}$ with
isolated singularities the number of exponents in
$(\alpha,\alpha+1]$ (counted with multiplicity) is
upper-semicontinuous for any $\alpha\in\bR$.
This gives a necessary condition for adjacent relation
of isolated hypersurface singularities, and implies a
counterexample to some conjecture about the adjacent relation.
(For a lower weight deformation of a weighted
homogeneous polynomial, this is due to Varchenko [50].)

\medskip\noindent
{\bf 2.6.~Steenbrink's conjecture} [47].
If $\dim\,$Sing$\,f=1$, and $g$ is generic with $dg\ne 0$, then
we have for $r\gg 0$
$$
\aligned
&\Sp(f+g^r,x)-\Sp(f,x)\\
\quad&=\hbox{$\sum_{k,j}t^{\alpha_{k,j}+
(\beta_{k,j}/m_kr)}(1-t)/(1-t^{1/m_kr})$,}
\endaligned
\leqno(2.6.1)
$$
where $m_k=\,$mult$_xZ_k$ with $Z_k$ the irreducible components of
(Sing$\,f)_{\red}$, the $\alpha_{k,j}$ are the exponents
(counted with multiplicities) at $y\in Z_k\setminus\{x\}$, and
$\beta_{k,j}$ are rational numbers in $(0,1]$ such that
$\exp(-2\pi i\beta_{k,j})$ are the eigenvalues of the monodromy
along $Z_k\setminus\{x\}$ (compatible with $\alpha_{k,j}$),
see [36] for a proof.

The formula (2.6.1) can be used for the calculation of $\Sp(f+g^r,x)$,
see [47].

\bigskip\bigskip
\centerline{\bf 3. Multiplier ideals and an extension theorem}

\bigskip\noindent
In this section we recall the definition of multiplier ideals
in the general case and give an extension theorem generalizing
\Mustata's formula in the case of hyperplane arrangements.

\medskip\noindent
{\bf 3.1.~Definition.}
Let $Z$ be a subvariety of a complex manifold or a smooth complex
algebraic variety $X$.
The multiplier ideal $\cJ(X,\alpha Z)$ for $\alpha\in\bQ_{>0}$ is
defined by
$$
g\in\cJ(X,\alpha Z)\Leftrightarrow |g|^2/(\sum |f_i|^2)^{\alpha}
\,\,\hbox{is locally integrable},
\leqno(3.1.1)
$$
where $f_1,\dots,f_r$ are local generators of the ideal of $Z$, or
$$
\cJ(X,\alpha Z)=\rho_*\omega_{\tX/X}(-\msum_i\,[\alpha m_i]
\tD_i)),
\leqno(3.1.2)
$$
where $\rho:(\tX,\tD)\to(X,Z)$ is an embedded resolution
such that $\rho^{-1}\cI_Z$ generates the ideal $\cI_{\tD}$ of
$\tD=\sum_i m_i\tD_i$.

Define for any $\alpha$ (with $0<\varepsilon\ll 1$)
$$
\cG(X,\alpha Z)=\cJ(X,(\alpha-\varepsilon)Z)/\cJ(X,\alpha Z).
\leqno(3.1.3)
$$
We say that $\alpha$ is a jumping number of $Z$ if and only if
$\cG(X,\alpha Z)\ne 0$.
Set
$$
\JN(Z)=\{\hbox{Jumping numbers of $Z$}\}\subset\bQ_{>0}.
\leqno(3.1.4)
$$

\medskip\noindent
{\bf 3.2.~Extension of multiplier ideals.}
Assume $X=Y\times\bC^r$ and $D=f^{-1}(0)$ with $\lambda^*f=f$
for $\lambda\in\bC^*$,
where the action of $\lambda$ is defined by
$$
\lambda\cdot(y,z_1,\dots,z_r)=(y,\lambda^{w_1}z_1,\dots,
\lambda^{w_1}z_r)\in Y\times\bC^r,
$$
with $w_i>0$.
For $y\in Y=Y\times\{0\}\subset X$, let
$$
\aligned
G^{>\alpha}\cO_{X,y}&=\{g\in\cO_{X,y}\,|\,v(g)>\alpha\}\,\,\,
\hbox{with}\\
v(\msum a_{\nu}\,z^{\nu})&=\min\{\msum_i \,w_i(\nu_i+1)\,|\,
a_{\nu}\ne 0\}.
\endaligned
$$
Let $X'=X\setminus (Y\times\{0\})$, $D'=D\cap X'$ with the inclusion
$j:X'\to X$. Then

\medskip\noindent
{\bf 3.3.~Theorem} [39]. {\it We have}
$$\cJ(X,\alpha D)_y=(j_*\cJ(X',\alpha D'))_y\cap
G^{>\alpha}\cO_{X,y}.$$

\medskip
This implies the following generalization of \Mustata's formula~[29]
in the case of hyperplane arrangements (see (5.17) below).

\medskip\noindent
{\bf 3.4.~Corollary.} {\it
Assume $D$ is the affine cone of a divisor $Z$ on $\bP^{n-1}$.
Let $d=\deg Z=\deg f$. Then
$$
\hbox{$\cJ(X,\alpha D)=I_0^k$ with $k=[d\alpha]-n+1$ if
$\alpha<\alpha'_{f,0}$},
\leqno(3.4.1)
$$
where $I_0$ is the ideal of $0$ and
$\alpha'_{f,0}=\min_{x\ne 0}\{\alpha_{f,x}\}$.
}

\medskip\noindent
{\bf 3.5.~Corollary.} {\it With the above assumption
$$
\hbox{$\dim F^{n-1}H^{n-1}(F_0,\bC)_{{\bf e}(-k/d)}=
\binom{k-1}{n-1}$ for $0<\frac{k}{d}<\alpha'_{f,0}$},
$$
and the same holds with $F$ replaced by $P$.
}

\medskip\noindent
{\bf 3.6.~Corollary.} {\it
With the above assumption, we have}
$$\alpha_{f} =\min\Big(\alpha'_{f,0},\frac{n}{d}\Big).$$

\bigskip\bigskip
\centerline{\bf 4. Relations in the hypersurface case}

\bigskip\noindent
In this section we explain certain relations among $b$-function,
spectrum and multiplier ideals in the hypersurface case.

\medskip\noindent
{\bf 4.1.~Theorem} (Budur [7]) {\it
Assume $\Sing\,f=\{x\}$. Then
}
$$
\aligned
&n_{f,\alpha}=\dim\cG(X,\alpha D)_x\,\,\,(\alpha\in(0,1)),\\
&\JN(D)\cap(0,1)=E_f\cap(0,1).
\endaligned
\leqno(4.1.1)
$$

\medskip\noindent
(This is generalized to the non-hypersurface case in Th.~(7.4).)

\medskip\noindent
{\bf 4.2.~Theorem} (Budur, S.\ [10]). {\it
Let $V$ denote also the induced filtration on
$\cO_X\subset \cO_X[\rd_t]\delta(f-t)$.
If $\alpha$ is not a jumping number,
$$
\cJ(X,\alpha D)=V^{\alpha}\cO_X.
\leqno(4.2.1)
$$
For $\alpha$ general we have for $0<\varepsilon\ll 1$
}
$$
\cJ(X,\alpha D)=V^{\alpha+\varepsilon}\cO_X,\quad
V^{\alpha}\cO_X=\cJ(X,(\alpha-\varepsilon) D).
\leqno(4.2.2)
$$

\medskip
Note that $V$ is left-continuous and
$\cJ(X,\alpha D)$ is right-continuous, i.e.
$$
V^{\alpha}\cO_X=V^{\alpha-\varepsilon}\cO_X,\quad
\cJ(X,\alpha D)=\cJ(X,(\alpha+\varepsilon) D).
\leqno(4.2.3)
$$

The proof of (4.2) uses the theory of bifiltered direct images [34], [35]
to reduce the assertion to the normal crossing case.

As a corollary we get another proof of the results of Ein, Lazarsfeld,
Smith and Varolin [18], and of Lichtin, Yano and Koll\'ar [26].

\bigskip
\noindent
{\bf 4.3.~Corollary.}
(i) $\JN(D)\cap(0,1]\subset R_f$ ({\it see} [18]).

\noindent
(ii) $\alpha_f= \min\JN(D)$ ({\it see} [26]).

\medskip
Setting $\alpha'_{f,x}=\min_{y\ne x}\{\alpha_{f,y}\}$,
we have a partial converse of Cor.~(4.3)(i) as follows.

\medskip\noindent
{\bf 4.4.~Theorem.} {\it
If $\xi f =f$ for a vector field $\xi$, then
}
$$
R_f\cap(0,\alpha'_{f,x})=\JN(D)\cap
(0,\alpha'_{f,x}).
\leqno(4.4.1)
$$

\medskip
This does not hold without the assumption on
$\xi$ nor without restricting to $(0,\alpha'_{f,x})$.

\medskip\noindent
{\bf 4.5.~Brieskorn lattice} (isolated singularities case).
The Brieskorn lattice [5] and its saturation are defined by
$$
\aligned
H''_f &= \Omega_{X,x}^n/df\wedge d\Omega_{X,x}^{n-2},\\
\tH''_f&=\msum_{i\ge 0}(t\rd_t)^iH''_f \subset H''_f[t^{-1}].
\endaligned
\leqno(4.5.1)
$$
These are finite $\bC\{t\}$-modules with a regular singular
connection.
Here $\Omega_X^{\ssbull}$ is analytic and $n=\dim X$.
Note that the action of $\partial_t^{-1}$ on $H''_f$ is well-defined
by $\partial_t^{-1}[\omega]=[df\wedge\xi]$ where
$\xi\in\Omega_{X,x}^{n-1}$ such that $d\xi=\omega$ in $\Omega_{X,x}^n$.

\medskip\noindent
{\bf 4.6.~Theorem} (Malgrange [27]). {\it
In the isolated singularity case, the reduced $b$-function
$\tb_f(s)$ coincides with the minimal polynomial of $-\rd_tt$ on
$\tH''_f/t\tH''_f$.
}

\medskip
(The above formula of Kashiwara on $b$-function (1.16.1) can be
proved by using this together with Brieskorn's calculation.)

\medskip\noindent
{\bf 4.7.~Asymptotic Hodge structures}
(Varchenko [49] and Scherk-Steenbrink [43]).
In the isolated singularity case, let $\cG_f$ be the
Gauss-Manin system $H''_f[\rd_t]$ (which is the localization of
$H''_f$ by the action of the microdifferential operator
$\rd_t^{-1}$).
Let $V$ be the filtration of Kashiwara and Malgrange on $\cG_f$.
Set $n=\dim X$.
Then
$$
\aligned
&F^pH^{n-1}(F_x,\bC)_{\lambda}=\Gr_V^{\alpha}H''_f\\
&\quad\text{for}\,\,\, p=[n-\alpha],\lambda=e^{-2\pi i\alpha},
\endaligned
\leqno(4.7.1)
$$
under the canonical isomorphism
$$
H^{n-1}(F_x,\bC)_{\lambda}=\Gr_V^{\alpha}\cG_f,
\leqno(4.7.2)
$$
together with $\rd_t^i:\Gr_V^{\alpha}\cG_f\simto\Gr_V^{\alpha-i}
\cG_f$ for $i\in\bZ$.
Note that Varchenko's filtration is defined by using
$t^{-i}$ instead of $\rd_t^i$.

The formula (4.7.1) can be generalized to the non-isolated
singularity case using mixed Hodge modules.

\medskip\noindent
{\bf 4.8.~Reformulation of Malgrange's formula.}
Set
$$
\aligned
&\tF^pH^{n-1}(F_x,\bC)_{\lambda}=\Gr_V^{\alpha}\tH''_f\\
&\quad\text{for}\,\,\, p=[n-\alpha],\lambda=e^{-2\pi i\alpha},
\endaligned
\leqno(4.8.1)
$$
under the canonical isomorphism (4.7.2).
Then
$$
\tm_{f,\alpha}=\hbox{deg(min poly}
(N\,|\,\Gr_{\tF}^pH^{n-1}(F_x,\bC)_{\lambda})),
\leqno(4.8.2)
$$
where min poly means the minimal polynomial.

\medskip\noindent
{\bf 4.9.~Remarks.} (i)
If $f$ has an isolated singularity, then,
as a corollary of the results of Malgrange [27], Varchenko [49],
Scherk-Steenbrink [43] explained in (4.6--7), we have
$$
\tR_f\subset \mcup_{0\le k<n}\,(E_f-k),\,\,\,
\tal_f=\min \tR_f=\min E_f.
\leqno(4.9.1)
$$

\medskip
(ii) If $f$ is weighted homogeneous with an isolated singularity,
then by the result of Kashiwara explained in (1.16) we have
$$
\tF=F,\quad \tR_f = E_f.
\leqno(4.9.2)
$$

\medskip
(iii) Let $g$ be a weighted homogenous polynomial with an isolated
singularity, and $h$ be a monomial $x^u$ with modified degree
$\beta:=\sum_i w_i u_i>1$ where $w_1,\dots,w_n$ are the weights
associated to $g$, i.e. $\sum_iw_ix_i\rd g/\rd x_i=g$.
Set $f=g+h$, and assume $h\notin(\partial g)$.
Then $E_f\ne \tR_f$.
Indeed, we have
$\partial_t[x^{u+v}dx]\in\tH''_f$ for any monomial $x^v$ since
$$
(\msum_i(v_i+1)w_i)\partial_t^{-1}[x^vdx]-t[x^vdx]=(\beta-1)
[x^{u+v}dx].
$$
We can apply this to a monomial $x^v$ such that $x^{u+v}$ generates
the highest modified degree part of $\bC[x]/(\partial g)$ which
is 1-dimensional.

\medskip\noindent
{\bf 4.10.~Example.} Let $f=x^5+y^4+x^3y^2$. Then
$$
\aligned
E_f&=\{\hbox{$\frac{i}{5}+\frac{j}{4}$}\,|\,1\le i \le 4,\,
1\le j\le 3\},\\
\tR_f&=\hbox{$E_f\cup\{\frac{11}{20}\}\setminus \{\frac{31}{20}$}\}.
\endaligned
$$
This is the simplest example such that
$E_f\ne \tR_f$.

\medskip\noindent
{\bf 4.11.~Relation with rational singularities} [37].
{\it Assume $D:=f^{-1}(0)$ is reduced.
Then $D$ has rational singularities if and only if $\tal_f>1$.
Moreover,
$$
\omega_D/\rho_*\omega_{\tD}\simeq F_{1-n}\varphi_f\cO_X,
$$
where $\rho:\tD\to D$ is a resolution of singularities.
}

\medskip
In the isolated singularities case, this was proved in [32]
using the coincidence of $\tal_f$ and the minimal exponent.

\medskip\noindent
{\bf 4.12.~Relation with the pole order filtration} [37].
Let $P$ be the pole order filtration on $\cO_X(*D)$, i.e.
$P_i=\cO_X((i+1)D)$ if $i\ge 0$, and $P_i=0$ if $i<0$.
Let $F$ be the Hodge filtration on $\cO_X(*D)$.
Then we have $F_i\subset P_i$ in general, and
$$
F_i=P_i\,\,\hbox{on a neighborhood of $x$ if
$\,i\le\tal_{f,x} - 1$}.
$$

\medskip\noindent
{\bf 4.13.~Remark.} In the case $X=\bP^n$, replacing $\tal_{f,x}$
with $[(n-r)/d]$ where $r = \dim\Sing \,D$ and $d = \deg D$,
the assertion was obtained by Deligne (unpublished).

\bigskip\bigskip
\centerline{\bf 5. Hyperplane arrangement case}

\bigskip\noindent
In this section we treat the case of hyperplane arrangements.

\medskip\noindent
{\bf 5.1.} Let
$D$ be a central hyperplane arrangement in $X=\bC^n$, i.e.
$D$ is an affine cone of a projective hyperplane arrangement
$Z\subset\bP^{n-1}$.
Let $f$ be the reduced equation of $D$ with $d=\deg f>n$.
Assume $D$ is not the pull-back of $D'\subset\bC^{n'}$ with
$ n'<n$. Then we have

\medskip\noindent
{\bf 5.2.~Theorem.} (i) $\max R_f<2-\frac{1}{d}$. (ii) $m_1=n$.

\medskip
For the proof of (i) we use a partial generalization of a solution of
Aomoto's conjecture due to Esnault, Schechtman, Viehweg, Terao,
Varchenko ([19], [42]) together with the following generalization
of Malgrange's formula in (4.8).

\medskip\noindent
{\bf 5.3.~Theorem} ({\rm Generalization of Malgrange's formula}) [39].
{\it There exists a pole order filtration $P$ on
$H^{n-1}(F_0,\bC)_{\lambda}$ satisfying the following property.

If $(\alpha+\bN)\cap R'_{f,0}=\emptyset$ with
$R'_{f,0}=\cup_{x\ne 0}R_{f,x}$, then
$$
\alpha\in R_f\,\,\,\iff\,\,\,
\Gr_P^pH^{n-1}(F_0,\bC)_{\lambda}\ne 0,
\leqno(5.3.1)
$$
where $p=[n-\alpha],\lambda=e^{-2\pi i\alpha}$.
}

\medskip
Using this, the proof of (5.2)(i) is reduced to
$$
P^iH^{n-1}(F_0,\bC)_{\lambda}=H^{n-1}(F_0,\bC)_{\lambda},
\leqno(5.3.2)
$$
for $i=n-1$ if $\lambda=1$ or $e^{2\pi i/d}$,
and $i=n-2$ otherwise.

\medskip\noindent
{\bf 5.4.~Construction of the pole order filtration P.}
Let $U=\bP^{n-1}\setminus Z$, and $F_0=f^{-1}(1)\subset\bC^n$.
Then $F_0$ is canonically identified with a cyclic $d$-fold covering
$\pi:\tU\to U$ ramified over $Z$.
Let $L^{(k)}$ be the local systems of rank 1 on $U$ such that
$\pi_*\bC=\moplus_{0\le i<d}L^{(k)}$
and $T$ acts on $L^{(k)}$ by $e^{-2\pi i k/d}$. Then we have
canonical isomorphisms
$$
H^j(U,L^{(k)})=H^j(F_0,\bC)_{{\bf e}(-k/d)},
\leqno(5.4.1)
$$
and $P$ is induced by the pole order filtration on
the meromorphic extension $\cL^{(k)}$ (see [11]) of
$L^{(k)}\otimes_{\bC}\cO_U$ over $\bP^{n-1}$, see [16], [39], [40].
This is closely related to [1] and also the following.

\medskip\noindent
{\bf 5.5.~Solution of Aomoto's conjecture} ([19], [42]). Let
$Z_i$ be the irreducible components of $Z\,\,(1\le i\le d)$.
Let $g_i$ be the defining equation of $Z_i$ on
$\bP^{n-1}\setminus Z_d$ for $i<d$, and set
$$
\omega=\msum_{i<d}\,\alpha_i\omega_i\,\,\,\hbox{with}\,\,\,
\omega_i=dg_i/g_i,\,\alpha_i\in\bC.
$$
Let $\nabla$ be the connection on $\cO_U$ defined by
$$
\nabla u=du+\omega\wdg u.
$$
Set $\alpha_d=-\sum_{i<d}\alpha_i$.
Then $H_{\rm DR}^{\ssbull}(U,(\cO_U,\nabla))$ is
calculated by the complex
$$
(\cA^{\ssbull}_{\alpha},\omega\wedge) \quad\text{with}\quad
\cA^p_{\alpha}=\msum \bC\omega_{i_1}\wdg\cdots\wdg\omega_{i_p},
\leqno(5.5.1)
$$
if the following condition is satisfied:
$$
\hbox{$\sum_{Z_i\supset L}\alpha_i\notin \bN\setminus\{0\}$ for any
{\it dense} edge $L\subset Z$},
\leqno(5.5.2)
$$
see (5.7) below for dense edges.

\medskip
For the proof of (5.2)(ii) we use

\medskip\noindent
{\bf 5.6.~Proposition.} {\it If
$\Gr^W_{2n-2}H^{n-1}(F_x,\bC)_{\lambda}\ne 0$, then we have
$N^{n-1}\psi_{f,\lambda}\bC\ne 0$.}

\medskip
(Indeed, by the definition of $W$, we have the isomorphism
$$
N^{n-1}:\Gr^W_{2n-2}\psi_{f,\lambda}\bC\simto
\Gr^W_0\psi_{f,\lambda}\bC,
$$
and the assumption of (5.6) implies
$\Gr^W_{2n-2}\psi_{f,\lambda}\bC\ne 0$.)

\medskip
Note that Proposition~(5.6) implies (5.2)(ii), since we have the
nonvanishing of $\omega_{i_1}\wdg\cdots\wdg\omega_{i_{n-1}}$ in
$$
\Gr^W_{2n-2}H^{n-1}(\bP^{n-1}\setminus Z,\bC)=
\Gr^W_{2n-2}H^{n-1}(F_x,\bC)_{1}.
$$

\medskip\noindent
{\bf 5.7.~Dense edges.}
Let $D=\cup_iD_i$ be the irreducible decomposition.
Then $L=\cap_{i\in I}D_i$ for $I\ne \emptyset$ is called an edge of
$D$.
An edge $L$ is called {\it dense} if
$\{D_i/L\,|\,D_i\supset L\}$ is indecomposable.
Here $\bC^n\supset D$ is decomposable if
$\bC^n=\bC^{n'}\times\bC^{n''}$ such that
$D$ is the union of the pull-backs from $\bC^{n'},
\bC^{n''}$ where $n',n''\ne 0$.

Let $m_L=\#\{D_i\,|\,D_i\supset L\}$, and for $\lambda\in\bC$
$$
\aligned
\cD\cE(D) &= \{\hbox{dense edges of }D\},\\
\cD\cE(D,\lambda)&=\{L\in\cD\cE(D)\,|\,\lambda^{m_L}=1\}.
\endaligned
$$
We say that $L,\,L'$ are {\it strongly adjacent} if
$L\subset L'$ or $L\supset L'$ or $L\cap L'$ is non-dense.
Let
$$
\aligned
m(\lambda)&=\max\{|S|\,\big|\,
S\subset\cD\cE(D,\lambda) \,\, \text{such that any two edges}
\\
&\qquad\text{belonging to $S$ are strongly adjacent}\}.
\endaligned
$$

\medskip\noindent
{\bf 5.8.~Theorem} [40]. {\it We have $m_{f,\alpha}\le m(\lambda)$
with $\lambda=e^{-2\pi i\alpha}$.}

\medskip\noindent
{\bf 5.9.~Corollary.}
{\it We have $R_f\subset\bigcup_{L\in\cD\cE(D)}\bZ m_L^{-1}$.}

\medskip\noindent
{\bf 5.10.~Corollary.} {\it Assume that {\rm GCD}$(m_L,m_{L'})=1$ for any
strongly adjacent $L,L'\in\cD\cE(D)$. Then $m_{f,\alpha}=1$ for any
$\alpha\in R_f\setminus \bZ$.}

\medskip
For Theorem~(5.2) we use the canonical embedded resolution of
singularities $\pi:(\tX,\tD)\to (\bP^{n-1},D)$, see [42].
This is obtained by blowing up along
the proper transforms of the dense edges.
Note that we have mult $\tD(\lambda)_{\rm red}\le m(\lambda)$,
where $\tD(\lambda)$ is the union of $\tD_i$ such that
$\lambda^{\tm_i}=1$ and $\tm_i=\hbox{mult}_{\tD_i}\tD$.

\medskip\noindent
{\bf 5.11.~Generic case.} If $D$ is a generic central
hyperplane arrangement, then
$$
b_f(s)=(s+1)^{n-1}\mprod_{j=n}^{2d-2}(s+\hbox{$\frac{j}{d}$})
\leqno(5.11.1)
$$
by U.~Walther [53] (except for the multiplicity of $-1$)
using a completely different method.

Note that Theorems (5.2) and (5.8) imply that the left-hand side
divides the right-hand side of (5.11.1), and the equality follows
using also (3.5).

\medskip\noindent
{\bf 5.12.~Explicit calculation.}
Let $\alpha=k/d$,
$\lambda=e^{-2\pi i\alpha}$ with $k\in \{1,\dots,d\}$.
If $\alpha\ge \alpha'_{f,0}:=\min_{x\ne 0}\{\alpha_{f,x}\}$,
we assume there is $I\subset\{1,\dots,d-1\}$
with $|I|=k-1$, and also the condition of [42] (i.e. (5.5.2) above)
is satisfied for
$$
\hbox{$\alpha_i=1-\alpha$ if $i\in I\cup\{d\}$, and $-\alpha$
otherwise}.
\leqno(5.12.1)
$$

Let $V(I)$ be the subspace of $H^{n-1}\cA^{\ssbull}_{\alpha}$
generated by
$$
\omega_{i_1}\wdg\cdots\wdg\omega_{i_{n-1}}\quad\text{for}\quad
\{i_1,\dots,i_{n-1}\}\subset I.
$$

\medskip\noindent
{\bf 5.13.~Theorem.} {\it With the above notation and assumptions,
we have for $\alpha=k/d$ and $\lambda=e^{-2\pi i\alpha}$ with
$k\in \{1,\dots,d\}$ the following.

\noindent
{\rm (a)} In case $k=d-1$ or $d$, we have
$\alpha\in R_f$, $\alpha+1\notin R_f$.

\noindent
{\rm (b)} In case $\alpha<\alpha'_{f,0}$, we have $\alpha\in R_f$
if and only if $k\ge d$.

\noindent
{\rm (c)} In case $\binom{k-1}{n-1}<\dim H^{n-1}(F_0,\bC)_{\lambda}$,
we have $\alpha+1\in R_f$.

\noindent
{\rm (d)} In case $\alpha<\alpha'_{f,0}$,
$\alpha\notin R'_{f,0}+\bZ$ and
$\binom{k-1}{n-1}=\chi(U)$, we have $\alpha+1\notin R_f$.

\noindent
{\rm (e)} In case $\alpha\ge \alpha'_{f,0}$ and $V(I)\ne 0$,
we have $\alpha\in R_f$.

\noindent
{\rm (f)} In case $\alpha\ge \alpha'_{f,0}$ and $V(I)=H^{n-1}
\cA^{\ssbull}_{\alpha}$, we have $\alpha+1\notin R_f$.
}

\medskip\noindent
{\bf 5.14.~Theorem} [40]. {\it 
Assume $n=3$, $\max\{\hbox{\rm mult}_zZ\,|\,z\in Z\}=3$, and
$d\le 7$. Let $\nu_3$ be the number of triple points of $Z$, and
assume $\nu_3\ne 0$. Then
$$
\hbox{$b_f(s)=(s+1)\prod_{i=2}^4(s+\frac{i}{3})\,
\prod_{j=3}^{r}(s+\frac{j}{d})$},
\leqno(5.14.1)
$$
with $r=2d-2$ or $2d-3$.
We have $r=2d-2$ if $\nu_3<d-3$. The converse holds for $d<7$.
In the case $d=7$, we have $r=2d-3$ if $\nu_3>4$.
However, $r$ can be both $2d-2$ and $2d-3$ if $\nu_3=4$.
}

\medskip\noindent
{\bf 5.15.~Remarks.} (i) We have $\nu_3 < d-3$ if and only if
we have
$$
\chi(U)=\hbox{$\frac{(d-2)(d-3)}{2}$}-\nu_3>
\hbox{$\frac{(d-3)(d-4)}{2}=\binom{d-3}{2}$}.
$$

\medskip\noindent
(ii) By (5.4.1) we have
$$
\chi(U)=h^2(F_0,\bC)_{\lambda}-h^1(F_0,\bC)_{\lambda}\,\,\,
\text{if}\,\,\,\lambda^d=1\,\,\,\text{with}\,\,\,\lambda\ne 1.
$$

\medskip\noindent
(iii) Let $\nu'_i$ be the number of $i$-ple points of
$Z':=Z\cap\bC^2$. Then we have by [6]
$$
\hbox{$b_0(U)=1,\quad b_1(U)=d-1,\quad b_2(U)=\nu'_2+2\nu'_3,$}
$$

\medskip\noindent
{\bf 5.16.~Example.} 
Assume $Z'$ is defined by $(x^2-y^2)(x^2-1)(y^2-1)=0$ in $\bC^2$
with $d=7$.
Then (5.14.1) holds with $r=11$, and $12/7\notin R_f$.
In this case we have
$$
\aligned
&b_1(U)=6, \,\,b_2(U)=9, \,\,\chi(U)=4,\\
&\hbox{$h^2(F_0,\bC)_{\lambda}=4$ if $\lambda^7=1\,\,\text{and}\,\,
\lambda\ne 1$.}
\endaligned
$$
Then $5/7\in R_f$ by (e) and $12/7\notin R_f$ by (f),
where $I^c$ corresponds to $(x+1)(y+1)=0$.
Note that $5/7$ is not a jumping number.

\medskip\noindent
{\bf 5.17.~Multiplier ideals of hyperplane arrangements.}
Let $m_L=\hbox{mult}_LD$, $r=\codim_XL$, and
$\cI_L$ be the ideal of an edge $L\subset X$.
Then by \Mustata~[29]
$$
\hbox{$\cJ(X,\alpha D)=\bigcap_L \cI_L^{[\alpha m_L]+1-r_L}$}.
\leqno(5.17.1)
$$
(This is generalized as in Cor.~(3.4) above.)

As for the spectrum, it does not seem easy to give
a combinatorial formula even for the generic case,
see e.g. [39], 5.6.

\bigskip\bigskip
\centerline{\bf 6. $b$-function of a subvariety}

\bigskip\noindent
In this section we define the $b$-function in the general
case and explain a relation with the multiplier ideals.

\medskip\noindent
{\bf 6.1.}
Let $Z$ be a closed subvariety of a complex manifold or
a smooth complex algebraic variety $X$.
Let $f=(f_1,\dots,f_r)$ be generators of the ideal of $Z$.
(We do not assume $Z$ reduced nor irreducible.)
Define the action of $ t_{j} $ on
$$
\cO_{X}\bigl[\hbox{$\frac{1}{f_1\cdots f_r}$}\bigr][s_{1},\dots, s_{r}]
\mprod_{i}f_{i}^{{s}_{i}},
$$
by
$ t_{j}(s_{i}) = s_{i} + 1 $ if
$ i = j $, and
$ t_{j}(s_{i}) = s_{i} $ otherwise.
Set
$$
s_{i,j} = s_{i}t_{i}^{-1}t_{j},\quad s=\msum_i \,s_i.
$$
Then $b_f(s)$ is the monic polynomial of the smallest degree such that
$$
b_{f}(s)\mprod_{i}f_{i}^{{s}_{i}} =
\msum_{k=1}^r P_{k}t_{k}\mprod_{i}f_{i}^{{s}_{i}},
\leqno(6.1.1)
$$
where
$ P_{k} $ belong to the ring generated by
$ \cD_{X} $ and
$ s_{i,j} $.

\medskip
Here we can replace $\prod_{i}f_{i}^{{s}_{i}}$ with
$\prod_{i}\delta(t_i-f_i)$,
using the direct image by the graph of $f:X\to\bC^r$.
Note that the existence of $b_f(s)$ follows from the theory
of the $V$-filtration of Kashiwara and Malgrange.
This $b$-function has appeared in work of Sabbah [31] and
Gyoja [20] for the study of $b$-functions of several variables.

\medskip\noindent
{\bf 6.2.~Theorem} (Budur, \Mustata, S.\ [8]). {\it
Let $c=\hbox{codim}_XZ$.
Then $b_Z(s):=b_f(s-c)$ depends only on $Z$, i.e. it is independent
of the choice of $X$, $f=(f_1,\dots,f_r)$, and also of $r$.
}

\medskip\noindent
{\bf 6.3.~Equivalent definition.}
The $b$-function
$ b_{f}(s) $ coincides with the monic polynomial of the smallest degree
such that
$$
b_{f}(s)\mprod_{i}f^{s_{i}} \in\msum_{|c|=1}
\cD_{X}[\bs]\,
\mprod_{c_i<0}\hbox{$\binom{s_{i}}{-c_{i}}$}
\mprod_{i}f_{i}^{s_{i}+c_{i}},
\leqno(6.3.1)
$$
where
$ c = (c_{1}, \dots, c_{r}) \in \bZ^{r} $ with
$ |c|:=\sum_{i}c_{i} = 1 $.
Here $\cD_{X}[\bs]=\cD_{X}[s_1,\cdots,s_r]$.

This is due to \Mustata, and is used in the monomial ideal case,
see (8.7) below.
Note that the well-definedness does not hold
without the term $\prod_{c_i<0}\binom{s_{i}}{-c_{i}}$.

\medskip
We denote also by $V$ the induced filtration under the inclusion
$$
\cO_X\subset i_{f+}\cO_X=\cO_X[\rd_1,\dots,\rd_r]\mprod_i
\delta(t_i-f_i).
$$

\medskip\noindent
{\bf 6.4.~Theorem} (Budur, \Mustata, S.\ [8]). {\it
For $\alpha\notin\JN(Z)$, we have
$$
\cJ(X,\alpha Z)=V^{\alpha}\cO_X.
\leqno(6.4.1)
$$
In general we have for any $\alpha$ \(with $0<\varepsilon\ll 1)$
}
$$
\cJ(X,\alpha Z)=V^{\alpha+\varepsilon}\cO_X,\quad
V^{\alpha}\cO_X=\cJ(X,(\alpha-\varepsilon) Z).
\leqno(6.4.2)
$$

\medskip\noindent
{\bf 6.5.~Corollary} (Budur, \Mustata, S.\ [8]). {\it
We have the inclusion}
$$
\JN(Z)\cap[\alpha_f,\alpha_f+1)\subset R_f.
\leqno(6.5.1)
$$

\medskip\noindent
{\bf 6.6.~Theorem} (Budur, \Mustata, S.\ [8]). {\it
Assume $Z$ is reduced and is a local complete intersection. Then
$Z$ has at most rational singularities if and only if $\alpha_f=r$ with multiplicity $1$.
}

\bigskip\bigskip
\centerline{\bf 7. Spectrum of a subvariety}

\bigskip\noindent
In this section we define the spectrum in the general case
and explain a relation with the multiplier ideals.

\medskip\noindent
{\bf 7.1.}
Let $Z$ be a closed subvariety of a complex manifold or a
smooth complex algebraic variety $X$, and $\cI_Z\subset\cO_X$ be
the ideal sheaf of $Z$.
The Verdier specialization [52] is defined by
$$
\Sp_Z\bQ_X=\psi_t\,\bR j_*\bQ_{X\times\bC^*},
\leqno(7.1.1)
$$
where
$$
j:X\times\bC^*\,(=\cSpec_X\cO_X[t,t^{-1}])\to\cSpec_X
(\moplus_{i\in\bZ}\,\cI_Z^{-i}\otimes t^i)
$$
is the inclusion to the total space of the
deformation to the normal cone
$$
N_ZX=\cSpec_Z(\moplus_{i\in\bN}\,\cI_Z^i/\cI_Z^{i+1}).
\leqno(7.1.2)
$$

Let $\Lambda$ be an irreducible component of the fiber
$(N_ZX)_z$ over $z\in Z$, and
$\xi\in\Lambda$ be a sufficiently general point of $\Lambda$
with the inclusion $i_{\xi}:\{\xi\}\to N_ZX$.
Set $c_{\Lambda}=\dim X-\dim\Lambda$.
Define the non-reduced spectrum and the (reduced) spectrum
$$
\aligned
\hSp(Z,\Lambda)&=\msum_{\alpha>0}\,n_{\Lambda,\alpha}
t^{\alpha},\\
\Sp(Z,\Lambda)&=\hSp(Z,\Lambda)-(-1)^{c_{\Lambda}}
t^{c_{\Lambda}+1},
\endaligned
$$
where
$$
\aligned
n_{\Lambda,\alpha}&=\msum_j\,(-1)^j\dim\Gr_F^pH^{j+c_{\Lambda}}
(i_{\xi}^*\Sp_Z\bC_X)_{\lambda}\quad\hbox{with}\\
p&=[c_{\Lambda}+1-\alpha],\,\,\lambda=\exp(-2\pi i\alpha),\\
\endaligned
\leqno(7.1.3)
$$

If $(N_ZX)_x$ is irreducible
(e.g. if $Z$ is a complete intersection), set
$$
\hbox{$\hSp(Z,x)=\hSp(Z,\Lambda)$, etc. for $\Lambda=(N_ZX)_x$.}
$$
This generalizes the definition for hypersurfaces.

\medskip\noindent
{\bf 7.2.~Remarks.} (i) In general, we have
$$
n_{\Lambda,\alpha}=0\,\,(\alpha\le 0),\,\,\,\,
n_{\Lambda,\beta}\ge 0\,\,(\beta\in(0,1]).
$$
In the isolated complete intersection singularity case, we have
$$
\tn_{x,\alpha}\ge 0\quad\hbox{with}\quad
\Sp(Z,x)=\msum_{\alpha}\,\tn_{x,\alpha}t^{\alpha},
$$
but symmetry and semicontinuity do not hold, see [17], [30], [48].

\medskip
(ii) In the isolated complete intersection singularity case,
our definition coincides with the one
by Ebeling and Steenbrink [17] except for $n_{x,\alpha}$ with
$\alpha\in\bZ$.
Indeed, they take generic 1-parameter smoothings
$$
\hbox{$f:X'\to\bC\,\,$ of $\,\,Z$,\quad $g:X''\to\bC\,\,$ of
$\,\,X'$,}
$$
and consider $\varphi_f\psi_g\bQ_{X''}[n]$ (where $n=\dim Z$)
together with the exact sequence
$$
0\to\tH^n(F_f,\bC)\to\varphi_f\psi_g \bQ_{X''}[n]\to
H^{n+1}(F_g,\bC)\to 0,
$$
where $\psi_g \bQ_{X''}|_{X'\setminus\{0\}}=\bQ$,
$(\psi_g \bQ_{X''})_0=\bR\Gamma(F_g,\bC)$.
The action of the monodromy on $H^{n+1}(F_g,\bC)$ is
associated to the functor $\varphi_f$, and is the identity.

\medskip\noindent
{\bf 7.3.} 
Let $\cI_Z$ be the ideal sheaf of $Z\subset X$.
For $z\in Z$ and $\beta\in(0,1]\cap\bQ$, let
$$
\aligned
\cM(\beta)=\moplus_{i\ge 0}\,\cG(X,(\beta+i)Z),\quad&
\ocA=\moplus_{j\ge 0}\,\cI_Z^j/\cI_Z^{j+1},
\\
\cM(\beta,z)=\cM(\beta)/\fm_{Z,z}\cM(\beta),\quad&
\ocA(z)=\ocA/\fm_{Z,z}\ocA.
\endaligned
\leqno(7.3.1)
$$
Then $\cM(\beta)$, $\cM(\beta,z)$ are graded modules over
$\ocA$, $\ocA(z)$, because
$$
(\cI_Z^j/\cI_Z^{j+1})\,\cG(X,(\beta+i)Z)\subset
\cG(X,(\beta+i+j)Z).
\leqno(7.3.2)
$$
For $z\in Z$ and an irreducible component $\Lambda$ of
$(N_ZX)_z=\Spec\,\ocA(z)$,
$$
\mu_{\Lambda,\beta}:=\dim_{\bC(\Lambda)}\cM(\beta,z)
\otimes_{\ocA(z)}\bC(\Lambda),
\leqno(7.3.3)
$$
where $\bC(\Lambda)$ is the function field of $\Lambda$.

\medskip\noindent
{\bf 7.4.~Theorem} (Dimca, Maisonobe, S. [14]). {\it
Let $\beta\in(0,1]\cap\bQ$.

\noindent
{\rm (i)} We have $0\le n_{\Lambda,\beta}\le\mu_{\Lambda,\beta}$
\(in particular, $n_{\Lambda,\beta}=0$ if $z\notin\Supp\,
\cM(\beta))$.

\noindent
{\rm (ii)} We have
$n_{\Lambda,\beta}=\mu_{\Lambda,\beta}$ if
$\Supp_{\ocA}\cM(\beta)\subset (N_ZX)_z$ on a neighborhood
of the generic point of $\Lambda$.
}

\medskip
(For hypersurfaces, this is due to Budur [7].)

\medskip\noindent
{\bf 7.5.~Corollary} (DMS [14]). {\it
If $n_{\Lambda,\alpha}\ne 0$ with $\alpha\in(0,1)$, then
there is a nonnegative integer $j_0$ such that
$\alpha+j\in\JN(Z)$ for any $j\ge j_0\in\bN$.
}

\medskip\noindent
{\bf 7.6.~Theorem} (DMS [14]). {\it
If $T$ is a transversal slice to a stratum
of a good Whitney stratification and $r=\codim\,T$, we have
}
$$
\hSp(Z,\Lambda)=(-t)^r\,\hSp(Z\cap T,\Lambda).
$$

\medskip
(For the constantness of the jumping numbers under a topologically
trivial deformation of divisors, see [15].)

\medskip\noindent
{\bf 7.7.~Remark.} Let
$E_{Z,\Lambda}=\{\alpha\,|\,n_{\Lambda,\alpha}\ne 0\}$. Then
$$
\hbox{$\bigcup_{\Lambda}\exp(-2\pi iE_{Z,\Lambda})\subset
\exp(-2\pi iR_{f,x})$},
$$
where $\Lambda$ runs over the irreducible components of $(N_ZX)_x$.
However, the equality does not always hold (e.g. if $f=x^2y$)
unless we take the union over the irreducible components
$\Lambda$ of $(N_ZX)_y$ for any $y\in Z$ sufficiently near $x$.

\bigskip\bigskip
\centerline{\bf 8. Monomial ideal case}

\bigskip\noindent
In this section we treat the monomial ideal case.

\medskip\noindent
{\bf 8.1.~Multiplier ideals.} Let
$ \fa \subset\bC[x] := \bC[x_{1},\dots, x_{n}] $ a monomial ideal.
We have the associated semigroup defined by
$$
\Gamma_{\fa} = \{u\in\bN^n\mid x^u\in\fa\}.
$$
Let
$ P_{\fa} $ be the convex hull of $\Gamma_{\fa}$ in
$ \bR_{\ge 0}^{n} $.
Set ${\mathbf 1}=(1,\dots,1)$, and
$$
U(\alpha):=\{\nu\in\bN^n\,|\,{\nu}+{\mathbf 1}\in
(\alpha+\varepsilon)P_{\fa}\,\,(0<\varepsilon\ll 1)\}.
$$
Let $Z$ be the subvariety of $X=\bC^n$ defined by $\fa$.

By Howald we have the following.

\medskip\noindent
{\bf 8.2.~Theorem}\,({\it Multiplier ideals})\,(Howald [21]).
{\it We have}
$$
\cJ(X,\alpha Z)=\msum_{\nu\in U(\alpha)}\, \bC\,x^{\nu}.
$$

\medskip\noindent
{\bf 8.3.~Corollary.} {\it
Set $\phi(\nu)=L_{\sigma}(\nu)$ if $\nu\in \hbox{\rm Cone}(0,\sigma)
:=\bigcup_{\lambda\ge 0}\lambda\sigma$,
where $L_{\sigma}$ is as in {\rm 8.4} below.
Then
$$
\JN(Z)=\{\phi(\nu)\,|\,\nu\in\bZ_{>0}^n\}.
$$
}

\medskip\noindent
{\bf 8.4.~Spectrum.}
For a maximal face $\sigma$ of $P_{\fa}$, set

\medskip
$L_\sigma$ : the linear function such that
$L_\sigma^{-1}(1)\supset\sigma$,

\medskip
$c_\sigma$ : the smallest positive integer such that
$c_\sigma L_\sigma\in\bZ[x]$,

\medskip
$e_\sigma=|G'_\sigma/G_\sigma|$,

\medskip\noindent
where $G'_\sigma=\bZ^{n}\cap L^{-1}_\sigma(0)$ and
$G_\sigma$ is generated by $\nu-\nu'$ with
$\nu,\nu'\in\Gamma_{\fa}\cap \sigma$.

\medskip\noindent
{\bf 8.5.~Theorem}\,({\it Spectrum})\,(Dimca, Maisonobe, S.~[14]).
{\it
We have a one-to-one correspondence between the maximal compact
faces $\sigma$ of $P_{\fa}$ and the irreducible components
$\Lambda$ of the fiber $(N_ZX)_0$, and
}
$$
\hSp(Z,\Lambda)=\msum_{i=1}^{c_\sigma}\,e_\sigma t^{i/c_\sigma}.
$$

\medskip\noindent
{\bf 8.6.~$b$-function.}
For a face $\sigma$ of $P_{\fa}$, set

\medskip
$V_\sigma$ : the linear subspace generated by $\sigma$,

\medskip
$M_\sigma$ : the subsemigroup generated by
$u-v$

\qquad
with $u\in\Gamma_{\fa}$, $v\in\Gamma_{\fa}\cap \sigma$,

\medskip
$M_\sigma'=v_0+M_\sigma$ with $v_0\in\Gamma_{\fa}\cap \sigma$
(independent of $v_0$),

\medskip
$R_\sigma=\{L_\sigma(u)\mid u\in ((M_\sigma\setminus M_\sigma')+
{\mathbf 1})\cap V_\sigma\}$,

\qquad
where ${\mathbf 1}=(1,\dots,1)$,

\medskip
$R_{\fa}=\{$roots of $b_{\fa}(-s)\}$ where $b_{\fa}(s)=b_Z(s)$.
 
\medskip\noindent
{\bf 8.7.~Theorem}\,({\it $b$-function})\,(Budur, \Mustata, S.\ [9]).
{\it We have $R_{\fa}=\bigcup_\sigma R_\sigma$ with $\sigma$
not contained in any coordinate hyperplanes.
}

\medskip\noindent
{\bf 8.8.~Remark.}
It is possible that $R_{\sigma}$ depends on the other $\sigma'$.
Indeed, we have the following (see [9], Ex.~4.4).

\medskip\noindent
(i) If $\fa=(xy^5,x^3y^2,x^5y)$, then
$R_{\fa}=R_\sigma\cup R_{\sigma'}$ with
$$
\hbox{$R_\sigma=\big\{\frac{5}{13},\frac{i}{13}\,
(7\le i\le 17),\frac{19}{13}\big\},\quad
R_{\sigma'}=\big\{\frac{j}{7}\,(3\le j\le 9)\big\}$.}
$$
So $R_\sigma=\big\{\frac{3i+2j}{13}\,(1\le i\le 3,1\le j\le 5)\big\}$
with $L_\sigma(i,j)=\frac{3i+2j}{13}$.

\medskip
(As for $R_{\sigma'}$ there is a misprint in loc.~cit.
as remarked by a student of W.~Veys.)

\medskip\noindent
(ii) If $\fa=(xy^5,x^3y^2,x^4y)$, then
$R_{\fa}=R_\sigma\cup R_{\sigma'}$ with
$$
\hbox{$R_\sigma=\big\{{\frac{i}{13}\,(5\le i\le 17)\big\},\quad
R_{\sigma'}=\big\{\frac{j}{5}\,(2\le j\le 6)\big\}}$.}
$$
So $R_\sigma\ne\big\{\frac{3i+2j}{13}\,(1\le i\le 3,1\le j\le 5)
\big\}$ with $\frac{19}{13}$ shifted to $\frac{6}{13}$.

\medskip\noindent
{\bf 8.9.~Comparison.}
Let $D=f^{-1}(0)$ for $f=\sum c_{\nu}x^{\nu}\in\bC[x]$
with non-degenerate Newton boundary $\rd P_f=\rd P_{\fa}$.
Assume $Z_{\red}=\{0\}$ so that $\Sing\,D=\{0\}$. Then
$$
\begin{array}{ccc}
\JN(D)\cap(0,1) & \overset{(1)}{=\!=} &
\JN(Z)\cap(0,1)\\
{\scriptstyle (2)}||\,\,\, & &
\cap{\scriptstyle (3)}\\
E_f\cap(0,1) & \underset{(4)}{\subset} &
\mcup_{\Lambda}E_{Z,\Lambda}\cap(0,1)\\
\end{array}
\leqno(4.5.3)
$$
where $E_{Z,\Lambda}=\{\alpha\,|\,n_{\Lambda,\alpha}\ne 0\}$.
Indeed, we have (1) by Howald [23], and (2) by Budur [7].
The composition of (1) and (2) is an equality by comparing the
formulas of Howald [21] and Steenbrink [45] (see also [33], [51]).
Finally we have (3) and (4) by [14].
(In general (3) (4) are not equality.)

\medskip\noindent
{\bf 8.10.~Example.}
If $\fa=(x_1^{m_1},\dots,x_n^{m_n})$, set
$$
\aligned
c_\sigma&=\hbox{LCM}(m_1,\dots,m_n),\quad
e_{\sigma}=m_1\cdots m_n/c_{\sigma},\\
E&=\{(a_1,\dots,a_n)\in\bN^n\,|\,a_i\in[1,m_i]\}.
\endaligned
$$
Then
$$
\aligned
\JN(Z)&=\Big\{\sum_{i=1}^n\frac{a_i}{m_i}
\,\Big|\,a_i\in\bZ_{>0}\Big\},\\
\hSp(Z,0)&=\sum_{i=1}^{c_{\sigma}}e_{\sigma}\,t^{i/c_{\sigma}},\\
b_{\fa}(s)&=\Big[\prod_{(a_1,\dots,a_n)\in E}\Big(s+\sum_{i=1}^n
\frac{a_i}{m_i}\Big)\Big]_{\red}.
\endaligned
$$
Here $[\prod_j(s+\beta_j)^{n_j}]_{\red}=\prod_j(s+\beta_j)$
if the $\beta_j$ are mutually different and $n_j\in\bZ_{>0}$.

\medskip
This may be compared with the following.

\medskip\noindent
{\bf 8.11.~Example.}
If $f=\sum_i x_i^{m_i}$ and $D=f^{-1}(0)$, set
$$
\tE=\{(a_1,\dots,a_n)\in\bN^n\,|\,a_i\in[1,m_i-1]\}.
$$
Then
$$
\aligned
\JN(D)\cap(0,1]&=\Big\{\sum_{i=1}^n\frac{a_i}{m_i}
\,\Big|\,a_i\in\bZ_{>0}\Big\}\cap(0,1],\\
\hbox{with}\,\,\,&\JN(D)=\big(\JN(D)\cap(0,1]\big)+\bN,\\
\Sp(D,0)&=\prod_{i=1}^n(t-t^{1/m_i})/(t^{1/m_i}-1),\\
\tb_f(s)&=\Big[\prod_{(a_1,\dots,a_n)\in \tE}\Big(s+\sum_{i=1}^n
\frac{a_i}{m_i}\Big)\Big]_{\red}.
\endaligned
$$
Indeed, for the assertion on $\JN(D)$, we can apply [22] or [7]
(i.e.\ Th.~(4.1) above), see also Th.~(4.4).
The other assertions follow from (1.16) and (2.3).
Note that the assertions hold for an isolated weighted homogeneous
singularities with weights $w_1,\dots,w_n$ if we replace $1/m_i$ by
$w_i$.

\medskip\noindent
{\bf 8.12.~Remark.}
In the monomial ideal case, $j_0$ in Cor.~(7.5) is bounded by
$n-1$, and $\JN(Z)$ is stable by adding any positive integers,
see [14].
Note that $j_0=n-1$ if the $m_i$ in (8.10) are mutually prime.
In general it is unclear whether $j_0$ is always bounded by $n-1$.

\bigskip\quad
\ver
\end{document}